\documentclass{sigcomm-alternate}
\usepackage{amsfonts,latexsym,amsmath,amstext,amssymb,verbatim,epsfig,psfrag}
\usepackage{colordvi,pstcol}
\usepackage{graphicx,psboxit}
\usepackage{psfrag}

\begin{document}
\title{D-iteration: application to differential equations}

\numberofauthors{1}
\author{
   \alignauthor Dohy Hong\\
   \affaddr{Alcatel-Lucent Bell Labs}\\
   \affaddr{Route de Villejust}\\
   \affaddr{91620 Nozay, France}\\
   \email{\normalsize dohy.hong@alcatel-lucent.com}
}

\date{\today}
\maketitle

\begin{abstract}
In this paper, we study how the D-iteration algorithm can be applied to numerically solve the differential equations such as heat equation in 2D or 3D. The method can be applied on the class of problems that can be addressed by the Gauss-Seidel iteration, based on the linear approximation of the differential equations.
\end{abstract}
\begin{psfrags}
%%%%%%%%%%%%%%%%%%%%%%%%%%%%%%%%%%%%%%%%%%%%%%%%%%%%%%%%%%%%
\section{Introduction}\label{sec:intro}
%%%%%%%%%%%%%%%%%%%%%%%%%%%%%%%%%%%%%%%%%%%%%%%%%%%%%%%%%%%%
The iterative methods to solve differential equations based on the linear approximation are very well
studied approaches \cite{Johnson_1987}, \cite{Ascher:1998:CMO:551054}, \cite{Podlubny:395913}, \cite{Gear:1971:NIV:540426},
\cite{Smith_1985}, \cite{Golub1996}, \cite{Saad}.
The approach we propose here (D-iteration) is a new approach initially applied to numerically solve the eigenvector
of the Pagerank type equation \cite{dohy}, \cite{d-algo}, \cite{dist-test}, \cite{distributed}, \cite{partition}, 
\cite{revisit}, \cite{revisit2}.

The D-iteration, as diffusion based iteration, is an iteration method that can be understood as a column-vector
based iteration as opposed to a row-vector based approach. Jacobi and Gauss-Seidel iterations are good examples of
row-vector based iteration schemes. While our approach can be associated to the {\em diffusion} vision, the
existing ones can be associated to the {\em collection} vision.

In this paper, we are interested in the numerical solution for linear equation:
\begin{eqnarray}\label{eq:le}
A.X &=& B
\end{eqnarray}
where $A$ and $B$ are the matrix and vector associated to the linear approximation of
differential equations with initial conditions or boundary conditions. 

While it is quite clear why diffusion approach is interesting (\cite{d-algo, revisit, revisit2}) when we consider
a sparse matrix with a very variable structure of in-degree/out-degree links (non-zero entries of the
matrix $A$), the problem statement in the context of linear system associated to differential equation 
is very different and we try to analyse/explain if there may be an interest to consider a solution
such as the D-iteration for those very regular sparse matrix.

%%%%%%%%%%%%%%%%%%%%%%%%%%%%%%%%%%%%%%%%%%%%%%%%%%%%%%%%%%%%
\section{Example of heat equation}\label{sec:heat}
%%%%%%%%%%%%%%%%%%%%%%%%%%%%%%%%%%%%%%%%%%%%%%%%%%%%%%%%%%%%

A typical linearized equation of the stationary heat equation in 2D is of the form:
\begin{align*}
T(n,m) &= \frac{1}{4} \left( T(n-1,m)+T(n,m-1)\right.\\
&\qquad \left.+T(n+1,m)+T(n,m+1)\right)
\end{align*}
which can be obtained by the discretization of the basic equation:
\begin{eqnarray}
\Delta.T(x,y) &=& \frac{\partial^2 T}{\partial x^2} + \frac{\partial^2 T}{\partial y^2} = 0.
\end{eqnarray}
inside the surface $\Omega$ (for instance, $\Omega = [0,L_x]\times[0,L_y]$)
Then additive terms appear for the initial or boundary conditions (Dirichlet) on the
frontier $\partial\Omega$ (for instance, for $x=0$ or $y=0$ etc).

For this family of equations, the linear dependences remain local (such as averaging of
neighbour positions' values).
From the intuitive point of view, the application of iterative methods on such a system
is convergent, because the system is diagonally dominant and strictly dominant for the 
boundary positions, so that the spectral radius of the global system is strictly less
than 1. As a consequence, compared to PageRank type equation with a damping factor $d$
for which the spectral radius of the system is explicitly equal to $d<1$ (per row), the
convergence of the iterative scheme should be slower for the same size of vector $N$.
Indeed, for the linear system we consider here, from the fluid diffusion point of view,
instead of having (PageRank equation) a fluid decreasing factor $1-d$ per 
entry level or vector level iteration,
the fluid can only disappear when it reaches the boundary positions (for instance, the
positions where the temperature is imposed from a heat source).

Now, in such a context, may the diffusion based iteration scheme be useful, faster?
We don't pretend to give the answer in this paper, we'll just illustrate the differences
and possible advantages through very simple examples.

The D-iteration requires updating two vectors: the fluid vector $F$ and the history
vector $H$ instead of a single vector for the Gauss-Seidel. It has been explained that
the $H$ vector is the exact connection to the Gauss-Seidel iterations (cf. \cite{revisit})
when starting from {\em empty} initial condition $H_0 = [0,..,0]$.
As it was the case for the iteration on the web graph associated matrix, the utility
of $F$ is to give us the exact information on the quantity of fluid that is sent.

\begin{figure}[htbp]
\centering
\includegraphics[width=6cm]{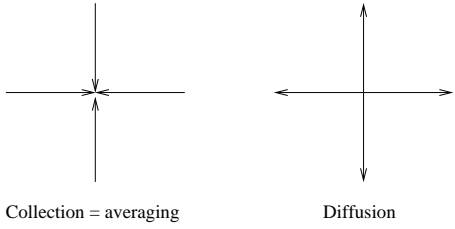}
\caption{Intuition: collection vs diffusion.}
\label{fig:cd}
\end{figure}

Whereas with the web graph, we may wish to apply the diffusion to the node having
the maximum amount of fluid (because the fluid that disappears is directly
proportional to this quantity), with the matrix we have here, we may wish to find
a way to push fluid to the boundary positions as quickly as possible, because inside
$\Omega$, the fluid does not disappear.

%Note that here, because of the spacial symmetry, the

%%%%%%%%%%%%%%%%%%%%%%%%%%%%%%%%%%%%%%%%%%%%%%%%%%%%%%%%%%%%
\section{First tests and algorithm adaptation}\label{sec:compa}
%%%%%%%%%%%%%%%%%%%%%%%%%%%%%%%%%%%%%%%%%%%%%%%%%%%%%%%%%%%%
As a first very simple illustration, we consider the iteration scheme:
\begin{align}
T(n,m) &= \frac{1}{4} \left( T(n-1,m)+T(n,m-1)\right.\\
&\qquad \left.+T(n+1,m)+T(n,m+1)\right)\label{eq:gs}
\end{align}
for $n= 1,..,L_x$ and $m= 1,.., L_y$ with the boundary condition:
\begin{align*}
T(n,m) &= 0, \mbox{ if } n=0 \mbox{ or } n=L_x\\
T(n,m) &= 0, \mbox{ if } m=L_y\\
T(n,m) &= 100, \mbox{ if } m=0.
\end{align*}

We first implemented the Gauss-Seidel iteration using the collection
operations defined by:
\begin{verbatim}
Collection of (n,m): (averaging)
 T[n][m] := 
    0.25*(T[n-1][m]+T[n][m-1]+T[n+1][m]+T[n][m+1]);
\end{verbatim}
when $(n,m)\notin \partial\Omega$
and the stopping condition: 
$$
\max_{n,m} |T(n,m)^k-T(n,m)^{k-1}| < \epsilon.
$$
Then, compared its computation cost to the D-iteration approach with the diffusion condition:
\begin{verbatim}
Diffusion of (n,m), if:
 F[n][m] > error;
\end{verbatim}

Remark that if there is no position satisfying $F[n][m] > error$, 
the condition becomes equivalent to the above stopping condition for Gauss-Seidel iteration.

For $L_x=L_y=1000$ and an error of $\epsilon = 0.1$ (let's say, we want to know the
temperature with a precision of 0.1 degree; one should be careful that the stopping
condition defined on the norm $\max$ $L_{\max}$ does not guarantee that the obtained result
is effectively at distance $\epsilon$ to the limit), we found the results:
\begin{itemize}
\item Gauss-Seidel (GS):
  2.35057e+08 operations of collections (over 236 cycles) in 8.3 seconds;
\item D-iteration (DI): 
  2.61681e+06 operations of diffusions (over 629 cycles) in 5.4 seconds.
\end{itemize}
\begin{verbatim}
Diffusion of (n,m):
 received := F[n][m]/4;
 T[n-1][m] += received;
 T[n][m-1] += received;
 T[n+1][m] += received;
 T[n][m+1] += received;
\end{verbatim}

In this case, the boundary condition is such that
the heat diffusion is progressive toward the position $y=L_y$.
Observe that with D-iteration, we did much less operations (by factor 100!), but this
implied the diffusion condition tests done on more cycles (629) than for Gauss-Seidel (236 cycles).

To reduce the diffusion condition test cost, we modified the diffusion condition (DC) as:
\begin{verbatim}
DC:
Diffusion of (n,m), if:
 F[n][m] > error/delta;
\end{verbatim}

The results when varying $\delta$ is shown in Table \ref{tab:1}:
we see a clear deterministic effect of $\delta$ when $\delta$ is a multiple
of 4 (due to the factor 1/4 of the matrix). 4 seems to be here the optimal value.

\begin{table}
\begin{center}
\begin{tabular}{|l|ccc|}
\hline
$\delta$        & cycles & diffusions   & time\\
\hline
1   & 629 & 2.61681e+06 & 5.4\\
2   & 317 & 6.65764e+06 & 3.2\\
3   & 232 & 6.07378e+06 & 6.1\\
4   & 234 & 6.56046e+06 & 2.5\\
5   & 235 & 6.89453e+06 & 6.2\\
6   & 235 & 7.11393e+06 & 6.2\\
7   & 235 & 7.30539e+06 & 6.2\\
8   & 235 & 7.45865e+06 & 2.5\\
16  & 236 & 8.2354e+06  & 2.5\\
\hline
\end{tabular}\caption{Impact of $\delta$ on the computation cost.}\label{tab:1}
\end{center}
\end{table}

If our understanding is correct, we should have more gain when $L_y$ is increased
(more stable value): as expected, this is confirmed in Table \ref{tab:2}.

\begin{table}
\begin{center}
\begin{tabular}{|l|cc|cc|}
\hline
$L_y$  & GS  &      &   DI   & \\
       & collections  &  time    & diffusions   & time\\
\hline
1000   & 2.3e+08 & 8.25  & 6.5e+06 & 2.47\\
2000   & 4.7e+08 & 12.68 & 6.5e+06 & 4.79\\
3000   & 7.0e+08 & 15.68 & 6.5e+06 & 7.14\\
4000   & 9.4e+08 & 18.76 & 6.5e+06 & 9.42\\
5000   & 1.1e+09 & 21.71 & 6.5e+06 & 11.77\\
\hline
\end{tabular}\caption{Impact of $L_y$ on the computation cost.}\label{tab:2}
\end{center}
\end{table}

While the  number of diffusions stays constant for DI, for GS the number of
collection operations increases proportionally to $L_y$. The increase of the
run time with DI is only due to the diffusion condition to be tested on a larger
set of points.

With this DC modification, we select the fluid diffusion only
at positions where the variation is worth to apply diffusion: when diffusion is
applied at position $(n,m)$, $F[n][m]$ is exactly the value by which
$H[n][m]$ is increased.
With the averaging computation, we are likely to do a lot of useless computations
in case certain points of space have converged under the desired precision.
And testing this condition would imply the same cost than applying the collection
operation
(because this would require the knowledge of the variation of the values of 
neighbour positions):
this condition test is exactly what we may do efficiently with $F$.

Here, we observe that the cost of testing diffusion condition (DC) is very important compared
to the diffusion cost (see the collections/diffusions ratio compared to the runtime ratio).
This is the main difference with the web graph: with web graphs (iterators required for
diffusions or collections), the dominant computation cost was the access time to the
entries of the matrix with the iterators. Here, this is no more the dominant cost.
Therefore, we need a specific optimization to reduce the DC tests.

In order to reduce the diffusion condition test cost, we introduced
a new (bool) state variable {\em open}:
\begin{verbatim}
bool open[n_x][n_y];
\end{verbatim}
which is set to true initially. 
Then, this state variable is updated as follows:
\begin{itemize}
\item when (DC) is not true, we set its {\em open} state to {\em false};
\item when an open position sends fluids to its neighbour positions, 
 these neighbour positions are set to {\em true}.
\end{itemize}

Then, we test this state value before testing the diffusion condition DC.
The results obtained by this modification is illustrate in Table \ref{tab:3}.
\begin{table}
\begin{center}
\begin{tabular}{|l|cc|cc|}
\hline
$L_y$  & GS  &      &   DI   & \\
       & collections   & time    & diffusions   & time\\
\hline
25     & 4.3e+06 & 0.07  & 4.0e+06 & 0.11\\
no op  &         & 0.3   &         & 0.44\\
\hline
50     & 1.1e+07 & 0.16  & 6.6e+06 & 0.2\\
no op  &         & 0.74  &         & 0.76\\
\hline
100    & 2.3e+07 & 0.34  & 6.6e+06 & 0.27\\
no op  &         & 1.5   &         & 0.86\\
\hline
1000   & 2.3e+08 & 8.25  & 6.6e+06 & 1.23\\
no op  &         & 21    &         & 2.3\\
\hline
2000   & 4.7e+08 & 12.68 & 6.6e+06 & 2.23\\
no op  &         & 37    &         & 3.9\\
\hline
3000   & 7.0e+08 & 15.68 & 6.6e+06 & 3.22\\
no op  &         & 51    &         & 5.5\\
\hline
4000   & 9.4e+08 & 18.76 & 6.6e+06 & 4.23\\
no op  &         & 66    &         & 7.0\\
\hline
5000   & 1.1e+09 & 21.71 & 6.6e+06 & 5.21\\
no op  &         & 82    &         & 8.6\\
\hline
\end{tabular}\caption{Speeding-up the DC test. GS: 236 cycles, DI: 234 cycles (except for $L_y=25$: 193, 193). no op: without compiler optimization (option $-O2$ with $g++$).}\label{tab:3}
\end{center}
\end{table}

We see that there is a significant gain on the runtime.
Finally, with the adaptation of DC and the introduction of {\em open} state
variables, the D-iteration may be an interesting candidate to efficiently
solve numerically the differential equations.
%Figure \ref{fig:} shows the comparison of the convergence w.r.t. the time.

%%%%%%%%%%%%%%%%%%%%%%%%%%%%%%%%%%%%%%%%%%%%%%%%%%%%%%%%%%%%
\section{Cost analysis}\label{sec:anal}
%%%%%%%%%%%%%%%%%%%%%%%%%%%%%%%%%%%%%%%%%%%%%%%%%%%%%%%%%%%%
\subsection{Cost decomposition}
In order to further understand the performance improvement and the
computation cost structure, we introduce the following assumption:
we assume that the runtime $RT$ of the solution computation is of the
form:
$$
RT = (\alpha + \beta)\times (L_x\times L_y)\times \mbox{nb\_iter} + c,
$$
where $\alpha$ is the unitary computation cost of one cycle iteration,
$\beta$ the average computation cost for one diffusion or one collection,
$c$ a constant cost and nb\_iter the number of iterations of the cycles.

To estimate the $\alpha$, $\beta$, $c$ values, we first iterated on $x$, $y$
without the collection operations:

\begin{verbatim}
sum = 0.0;
counter = 0;
while ( counter < nb_iter ){
  counter++;
  for (int x=0; x < L_x; x++){
    for (int y=0; y < L_y; y++){
      if ( !is_boundary(x,y) )
        sum += T[x][y];
    }
  }
}
\end{verbatim}
where \verb+is_boundary(x,y)+ is a variable (bool) that return true if we are
on the position with boundary conditions (or outside $\Omega$). 
Then $\alpha$ is estimated by dividing the runtime by nb\_iter and $L_x\times L_y$.
For $\beta$, we introduce collection operations and then $\alpha+\beta$ is estimated
by dividing the runtime by nb\_iter and $L_x\times L_y$.

Finally, $c$ is approximated by the initialization time.
For DI, we do the same evaluation with:
\begin{verbatim}
sum = 0.0;
counter = 0;
while ( counter < nb_iter ){
  counter++;
  for (int x=0; x < L_x; x++){
    for (int y=0; y < L_y; y++){
      if ( !is_boundary(x,y) and open[x][y] ){
        transit = F[x][y];
        if ( transit > error/scale ){
          sum += T[x][y];
        }
      }
    }
  }
}
\end{verbatim}
adding the DC test for $\alpha$ and then adding the diffusion operations
to estimate $\alpha+\beta$.
The results are shown on Table \ref{tab:cc}.

\begin{table}
\begin{center}
\begin{tabular}{|l|cc|cc|}
\hline
$L_x\times L_y$  &  GS    & &   DI &   \\
                 & $\alpha$  &  $\alpha+\beta$ & $\alpha$   & $\alpha+\beta$   \\
\hline
%$100\times 100$ & 4.3e-09 & 8.2e-09  & 8.0e-09 & 1.1e-08\\
%$200\times 200$ & 4.4e-09 & 8.6e-09 & 8.4e-09 & 1.3e-08\\
%$100\times 1000$   & 5.8e-09 & 8.8e-09 & 1.1e-08 & 1.4e-08\\
%$1000\times 1000$  & 4.9e-09 & 9.0e-09 & 8.9e-09 & 1.5e-08\\
%$200\times 5000$   & 4.5e-09 & 2.0e-08 & 1.2e-08 & 1.0-08\\
%$5000\times 200$   & 6.1e-09 & 8.7e-09 & 1.0e-08 & 1.4e-08\\
$100\times 100$ & 2.0e-08 & 6.2e-08  & 3.3e-08 & 1.0e-07\\
$200\times 200$ & 2.0e-08 & 6.1e-08  & 3.3e-08 & 1.0e-07\\
$1000\times 1000$  & 2.0e-08 & 6.7e-08  & 3.3e-08 & 1.0e-07\\
$100\times 1000$   & 2.0e-08 & 6.4e-08  & 3.3e-08 & 1.0e-07\\
$200\times 5000$   & 2.0e-08 & 6.8e-08  & 3.4e-08 & 1.0e-07\\
$5000\times 200$   & 2.2e-08 & 6.3e-08  & 3.4e-08 & 1.0e-07\\
\hline
\end{tabular}\caption{Computation cost. Without compiler optimization.}\label{tab:cc}
\end{center}
\end{table}

We see that $\alpha$ and $\beta$ are quite stable.
We see that the collection cost is roughly 2-3 times the iteration cost
for GS (we observed this factor is the one that varies the most with nb\_iter). 
For DI, its iteration cost is more than 1.5 times the
iteration cost for GS and the diffusion
cost is about 2 times the iteration cost (very stable with nb\_iter).

The runtime for GS in Table \ref{tab:3} for $1000\times 1000$ (no op) can
be decomposed as 5s (iterations) $+$ 15s (collections). In this case, DI converged
with about same number of cycles but with 30 times less diffusions (than collections).
So based on the above model, we would have a runtime of: 7.5s for iterations
and $15/30\times 2 = 1$s. But we observed 2.3s instead of 8.5s because most of time
we don't need to test twice the boundary and open condition. We observed that
only testing open condition is close (in this case) to the cost to test the
boundary condition, but also the main difference comes from the fact that when
the open condition is not true, we don't have to do any computation and
this had in this case an impact of cost reduction by up to 3-4 (the iteration
runtime for all values $open[x][y]$ set to false runs 3-4 times faster).
Therefore, the computation time obtained for DI would be: $5/4+1$.

Now, we could optimize a bit more the DC condition and the boundary position tests
including the boundary position information in the state variable {\em open}
(for instance using, 3 states variable: 0 for boundary position, 1 for not open to
DC test, 2 for open to DC test). Results are shown in Table \ref{tab:4}.

\begin{table}
\begin{center}
\begin{tabular}{|l|cc|cc|}
\hline
$L_y$  & GS  &      &   DI   & \\
       & collections   & time    & diffusions   & time\\
\hline
50     & 1.1e+07 & 0.16  & 6.6e+06 & 0.18\\
\hline
100    & 2.3e+07 & 0.34  & 6.6e+06 & 0.26\\
\hline
1000   & 2.3e+08 & 8.25  & 6.6e+06 & 1.13\\
\hline
2000   & 4.7e+08 & 12.68 & 6.6e+06 & 2.01\\
\hline
3000   & 7.0e+08 & 15.68 & 6.6e+06 & 2.89\\
\hline
4000   & 9.4e+08 & 18.76 & 6.6e+06 & 3.77\\
\hline
5000   & 1.1e+09 & 21.71 & 6.6e+06 & 4.68\\
\hline
\end{tabular}\caption{Speeding-up the DC test. GS: 236 cycles, DI: 234 cycles.}\label{tab:4}
\end{center}
\end{table}

\subsection{The number of diffusions per position}

Figure \ref{fig:position} shows the number of diffusions applied per position $y$.
We see that above $y=41$, the diffusion is no more applied, whereas with GS, the
collection operations need to be applied on all $y$ positions at each cycle.

\begin{figure}[htbp]
\centering
\includegraphics[angle=-90,width=\linewidth]{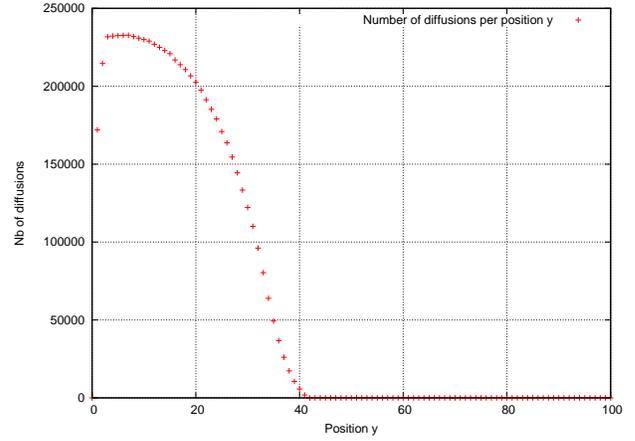}
\caption{Nb of diffusions applied: per y position.}
\label{fig:position}
\end{figure}

This example is of course for illustration of a situation when all positions does not
converge to the limit at the same speed.
If $L_y$ is closer to 41, there is no more gain using DC, because all positions are
converging almost at the same speed.

Figure \ref{fig:position2} shows the number of diffusions per position $y$
when in the case $L_y=1000$ we added 10 random positions where $T$ is imposed 
($T$ set to a random value between 0 and 1000): in this case, we obtained:
\begin{itemize}
\item GS: 697 cycles with 6.94e+08 collections for runtime of 9.5s;
\item DI: 687 cycles with 3.97e+07 diffusions for runtime of 4.0s.
\end{itemize}

\begin{figure}[htbp]
\centering
\includegraphics[angle=-90,width=\linewidth]{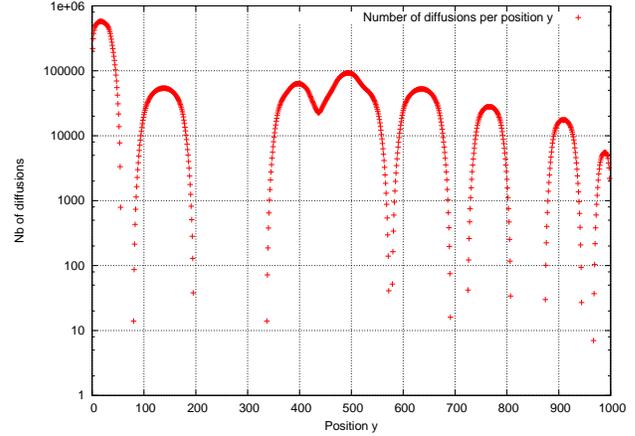}
\caption{Nb of diffusions applied: per y position: adding 10 random boundary conditions inside $\Omega$.}
\label{fig:position2}
\end{figure}

%Figure \ref{fig:position3} shows the results when the boundary condition is constant to 100.
%Which means that the limit must be $T=100$ constant everywhere.

%\begin{figure}[htbp]
%\centering
%\includegraphics[angle=-90,width=\linewidth]{nb_diffusion3.eps}
%\caption{Nb of diffusions applied: $T=100$ on the frontier of $\Omega$.}
%\label{fig:position3}
%\end{figure}

%In this case, we had:
%\begin{itemize}
%\item GS: 374 cycles with 3.72e+08 collections for runtime of 5.8s;
%\item DI: 372 cycles with 4.40e+07 diffusions for runtime of 2.9s.
%\end{itemize}

Figure \ref{fig:position4} is an illustration of the limit for
the boundary conditions of: $T=100$ on the frontier and $T=0$ imposed at 10 random
positions.

\begin{figure}[htbp]
\centering
\includegraphics[angle=-90,width=\linewidth]{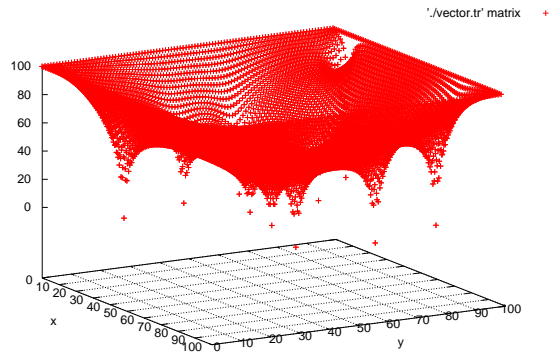}
\caption{Nb of diffusions applied: $T=100$ on the frontier of $\Omega$ and with 10 random 'holes' $T=0$.}
\label{fig:position4}
\end{figure}

%%%%%%%%%%%%%%%%%%%%%%%%%%%%%%%%%%%%%%%%%%%%%%%%%%%%%%%%%%%%
\section{Conclusion}\label{sec:conclusion}
%%%%%%%%%%%%%%%%%%%%%%%%%%%%%%%%%%%%%%%%%%%%%%%%%%%%%%%%%%%%
In this paper we addressed a first analysis of the potential of the
D-iteration when applied in the context of the numerical solving of
differential equations. We showed that this context requires an adaptation
of the D-iteration's fluid diffusion condition.
The results are quite promising and we hope to investigate further this
application case in the future.

\end{psfrags}
%======================================================================
\bibliographystyle{abbrv}
\bibliography{sigproc}

\end{document}